\documentclass[12pt]{article}

\newtheorem{theo}{Theorem}

\newtheorem{cor}{Corollary}

\def\R{{\bf R}}      \def\Z{{\bf Z}} \def\k{{\bf k}}
\def\z{\zeta}   \def\S{\zeta^{*}}
\def\={\;=\;} \def\-{\,-\,}
\def\I{I^{\text{adm}}}  \def\I{I_0}  \def\X{X_0}
\begin{document}
\title{Sum relations for multiple zeta values and connection formulas
for the Gauss hypergeometric functions}
\author{Takashi Aoki\footnote{Suppoted in part by JSPS Grant-in-Aid  
No.~14340042 and by No.~15540190}
\ \ and \ \ Yasuo Ohno\footnote{Suppoted in part by JSPS Grant-in-Aid   
No.~15740025 and by No.~15540190}}  \date{}
\maketitle

\begin{abstract} We give an explicit representation for the sums of
multiple zeta-star values of fixed weight and height in terms of
Riemann zeta values.\par
\end{abstract}

\section{Introduction}
In this article we establish a new family of relations between sums of
multiple zeta values and Riemann zeta values. This family
contains relations which do not appear in the family of
relations given in \cite{h-o},\cite{ohno1}, and \cite{o-z}.

Concerning multiple zeta values, there are two types of definition:  
multiple zeta values
  ``without equality'' and ones ``with equality'' (see below).
The former is mainly used in mathematical literature and the latter is  
the main subject
of this article. Normally multiple zeta values ({\it MZVs} for short)  
mean the former
and are denoted by $\zeta({\bf k})$.
We tentatively call the latter {\it multiple zeta-star values} and  
denote them by
$\zeta^*({\bf k})$ to distinguish
them from ordinary ones. We abbreviate them to {\it MZSVs}.
They are classic objects although there had been no name
of them. In fact, Euler was the first mathematician who was interested  
in
multiple zeta values and he mainly treated MZSVs \cite{euler}.
Recently, Hoffman \cite{hof} pointed out the significance of  
considering MZSVs
  as well as MZVs. The notation $S$ was used there for $\S$.

The main result of this article shows that the sum of MZSVs with fixed  
weight
and height turns out to be a rational multiple of Riemann zeta value at  
the weight.
Considering MZSVs clarifies the importance of those two indices:
weight and height. They have been played a role in \cite{l-m},  
\cite{o-z}.
The employment of the indices and MZSVs is a
neat way to formulate systematic description of relations that hold  
among
MZVs.
Another important index is depth.
We believe that MZSVs and the
three indices: weight, height and depth will play an important role in
investigation of
the structure of ${\bf Q}$-algebra generated by MZVs.
(Note that
this algebra coincides with  ${\bf Q}$-algebra generated by MZSVs.)

An interesting feature of the method employed in our proof is related  
to the theory
of differential equations in the complex domain. The method is a  
variation
on \cite{o-z} and the use of connection formulas of the Gauss  
hypergeometric
function is essential in both cases (see \cite{kombu}, \cite{ok-ue}  
also).

\section{Statement of the result}
For any multi-index $\k=(k_1,k_2,\ldots,k_n)$ ($k_i\in\Z$, $k_i>0$,
the {\it weight}, {\it depth}, and {\it height} of $\k$ are by  
definition
the integers $k=k_1+k_2+\cdots+k_n$, $n$, and $s=\#\{i|k_i>1\}$,  
respectively.
We denote by $I(k,s)$ the set of multi-indices $\k$ of weight $k$
and height $s$, and by $\I(k,s)$ the subset of {\it admissible}
indices, i.e., indices with the extra requirement that $k_1 \geq 2$.  
For any
admissible index $\k=(k_1,k_2,\ldots,k_n)\in\I(k,s)$, the {\it multiple
zeta values} $\S(\k)$ and $\z(\k)$ are defined by
$$ \S(\k)=\S(k_1,k_2,\ldots ,k_n) = \sum_{m_1\geq m_2 \geq \cdots \geq  
m_n \geq 1}
\frac{1}{{m_1}^{k_1} {m_2}^{k_2} \cdots {m_n}^{k_n}} \;.$$
$$ \z(\k)=\z(k_1,k_2,\ldots ,k_n) = \sum_{m_1>m_2 > \cdots > m_n>0}
\frac{1}{{m_1}^{k_1} {m_2}^{k_2} \cdots {m_n}^{k_n}} \;.$$
Note that, there are linear relations among $\S$ and $\z$, for example,
$$
\S(k_1,k_2)=\z(k_1,k_2)+\z(k_1+k_2),\ \  
\z(k_1,k_2)=\S(k_1,k_2)-\S(k_1+k_2),
$$
$$
\S(k_1,k_2,k_3)=\z(k_1,k_2,k_3)+\z(k_1+k_2,k_3)+\z(k_1,k_2+k_3)+\z(k_1+k 
_2+k_3),
$$
$$
\z(k_1,k_2,k_3)=\S(k_1,k_2,k_3)-\S(k_1+k_2,k_3)- 
\S(k_1,k_2+k_3)+\S(k_1+k_2+k_3),
$$
and so on. Multiple zeta-star values $\S$ had been studied by  
Euler\cite{euler}, and his study is
the origin of various researches of
multiple zeta values $\z$.
We consider a sum of multiple zeta-star values of fixed weight and  
height:
\begin{equation}\label{thesum}
\sum_{\k\in\I(k,s)}\S(\k)\;.
\end{equation}

Our main result will then be
\begin{theo} The sum {\rm (\ref{thesum})} is given by
\begin{equation}\label{maintheo}
\sum_{\k\in\I(k,s)}\S(\k)\= 2 {{k-1}\choose{2s-1}} (1-2^{1-k}) \z(k).
\end{equation}
\end{theo}

As an application of Theorem 1, we can express sums of
special values of Arakawa-Kaneko zeta function in terms of
Riemann zeta values.

For any positive integer $k \geq 1$,
T.~Arakawa and M.~Kaneko \cite{ak} defined the function $\xi_k(s)$ by
\[ \xi_k(s) = \frac{1}{\Gamma (s)}
\int_0^\infty \frac{t^{s-1}}{e^t-1} {\rm Li}_k(1-e^{-t})dt, \]
where ${\rm Li}_k(s)$ denotes the $k$-th polylogarithm
${\rm Li}_k(s)=\sum_{m=0}^\infty \frac{s^m}{m^k}$.
The integral converges for ${\rm Re}(s)>0$ and the function
$\xi_k(s)$ continues to an entire function of whole $s$-plane.
They proved that the special values of $\xi_k(s)$ at non-positive
integers are given by poly-Bernoulli numbers and the values at positive
integers
are given in terms of multiple zeta values.
Thereafter the second author \cite{ohno1} gave the following relation  
among
the values
of $\xi_k(s)$ at positive integers and MZSVs:
\[
\xi_k(n)=\S(k+1,\underbrace{1,\ldots ,1}_{\mbox{$n-1$}})
\]
where the both indices $k$ and $n$ are positive integers.
By using this relation, we have the following corollary of Theorem 1.

\begin{cor}
For any positive integer $k$, we have
\[
\sum_{n=1}^{k-1}\xi_{k-n}(n)=2(k-1)(1-2^{1-k})\z(k).
\]
\end{cor}

\section{Proof of Theorem 1}
We denote by $\X(k,s)$ the left-hand side of (\ref{maintheo}):
\begin{equation}\label{x}
\X(k,s)\= \sum_{\k\in\I(k,s)}\S(\k)\;.
\end{equation}
Since the set $\I(k,s)$ is non-empty only if the indices $k$ and $s$
satisfy the inequalities  $s\ge1$ and $k\ge 2s$, we can collect
all the numbers $\X(k,s)$ into a single generating function
\begin{equation}\label{phi0}
\Phi_0(x,z)\=\sum_{k,\,s}\X(k,s)\,x^{k-2s}\,z^{2s-2}
\quad\in\,\R[[x,z]]\,.
\end{equation}
Following \cite{o-z} and \cite{zag2}, we consider
  the multiple zeta-star value $\S(\k)$
  as the limiting value at $t\=1$ of the function
$$
L^*_{\k}(t)\= L^*_{k_1,k_2,\ldots ,k_n}(t) \=
\sum_{m_1 \geq m_2 \geq \cdots \geq m_n \geq 1}
\frac{t^{m_1}}{{m_1}^{k_1} {m_2}^{k_2} \cdots  
{m_n}^{k_n}}\qquad(|t|<1)\,.
$$
Note that we consider $L^*_{\k}(t)$ not just for $\k\in\I$ but for all  
$\k\in I$.
For $\k$ empty we define $L^*_{\k}(t)$ to be 1.
For non-negative integers $k$ and $s$ set
$$ X(k,s;t)\=\sum_{\k\in I(k,s)} L^*_{\k}(t)$$
(so $X(0,0;t)=1$ and $X(k,s;t)=0$ unless $k\geq 2s$ and $s\ge0$),
and let $\X(k,s;t)$ be the function defined by the same
formula but with the summation restricted to $\k\in\I(k,s)$.

We denote by
$\Phi=\Phi(x,z;t)$ and $\Phi_0=\Phi_0(x,z;t)$ the corresponding
generating functions
$$\Phi=\sum_{k,s \geq 0} X(k,s;t)x^{k-2s}z^{2s}
=1+L^*_1(t)x+L^*_{1,1}(t)x^2+\cdots$$
and
$$\Phi_0=\sum_{k,s \geq 0} \X(k,s;t)x^{k-2s}z^{2s-2}
=L^*_2(t)+L^*_{1,2}(t)x+L^*_3(t)x+\cdots\,.$$
Note that the coefficient of $x^{k-2s}z^{2s-2}$ in  
$\Phi_0(x,z;1)=\Phi_0(x,z)$ is $X_0(k,s)$.
Using the formula
$$ \frac{d}{dt}L^*_{k_1,\ldots,k_n}(t)\= \left\{
\begin{array}{ll}
{\displaystyle\frac{1}{\,t\,}}\,L^*_{k_1-1,k_2,\ldots,k_n}(t) &
  \quad{\rm if}\quad k_1\ge2,\\
  & \\
{\displaystyle\frac{1}{t(1-t)}}\,L^*_{k_2,k_3,\ldots,k_n}(t) &
  \quad{\rm if}\quad k_1=1 \\
\end{array}\right. $$
for the derivative of $L^*_{\bf k}(t)$, we obtain
$$\frac{d}{dt}\X(k,s;t)=
\frac 1t\Bigl(X(k-1,s-1;t)-\X(k-1,s-1;t)+\X(k-1,s;t)\Bigr)\,,$$
$$\frac{d}{dt}\Bigl(X(k,s;t)-\X(k,s;t)\Bigr)=\frac{1}{t(1-t)}X(k- 
1,s;t)\,,$$
or, in terms of generating functions,
$$\frac{d\Phi_0}{dt}=\frac1{xt}\Bigl(\Phi-1-z^2\Phi_0\Bigr)+\frac  
xt\Phi_0\,,
\qquad\frac d{dt}\Bigl(\Phi-z^2\Phi_0\Bigr)
=\frac{x}{t(1-t)}(\Phi-1)+\frac{x}{1-t}\,.$$
Eliminating $\Phi$, we obtain the differential equation
\begin{equation}\label{dfe}
t^2(1-t)\,\frac{d^2\Phi_0}{dt^2}+t\Bigl((1-t)(1-x)-x\Bigr)
\,\frac{d\Phi_0}{dt}+(x^2-z^2)\,\Phi_0=t
\end{equation}
for the power series $\Phi_0$.
The unique power-series solution at $t=0$ is given by
$$
\Phi_0(x,z;t)=\sum_{n=1}^\infty a_n t^n
$$
with
$$
a_n=\frac{\Gamma(n)\Gamma(n-x)\Gamma(1-x-z)\Gamma(1-x+z)}
{\Gamma(1-x)\Gamma(1-x-z+n)\Gamma(1-x+z+n)}.
$$
Here $\Gamma(z)$ denotes the gamma function.
Specializing to $t=1$ gives
\begin{equation}\label{phi01}
\Phi_0(x,z;1)=\sum_{n=1}^\infty a_n.
\end{equation}
We need to evaluate the right-hand side of (\ref{phi01}).
We can rewrite $a_n$ in the form
$$
a_n=\sum_{l=1}^{n}\left(  
\frac{A_{n,l}^{(+)}}{x+z-l}+\frac{A_{n,l}^{(-)}}{x-z-l} \right)
$$
with
$$
A_{n,l}^{(\pm)}=(-1)^l { n-1\choose l-1 }
\frac{(\pm z-l+1)(\pm z-l+2)\cdots(\pm z-l+n-1)}
{(\pm 2z-l+1)(\pm 2z-l+2)\cdots(\pm 2z-l+n)}.
$$
Hence we have
\begin{eqnarray*}
\sum_{n=1}^{\infty}a_n &=&
\sum_{n=1}^{\infty} \sum_{l=1}^{n}
\left( \frac{A_{n,l}^{(+)}}{x+z-l}+\frac{A_{n,l}^{(-)}}{x-z-l} \right)  
\\
&=& \sum_{l=1}^{n}
\left( \sum_{n=l}^{\infty} A_{n,l}^{(+)} \frac{1}{x+z-l}+
\sum_{n=l}^{\infty} A_{n,l}^{(-)} \frac{1}{x-z-l} \right) .
\end{eqnarray*}
The sums of $A_{n,l}^{(\pm)}$ in $n$ are evaluated as follows:
\begin{eqnarray*}
\sum_{n=l}^{\infty} A_{n,l}^{(\pm)}
&=& (-1)^l \sum_{n=0}^{\infty}
\frac{(l-1+n)!(\pm z-l+1)(\pm z-l+2)\cdots(\pm z+n-1)}{n!(l-1)!(\pm  
2z-l+1)(\pm 2z-l+2)\cdots(\pm 2z+n)} \\
&=& (-1)^l
\frac{(\pm z-l+1)(\pm z-l+2)\cdots(\pm z-1)}
{(\pm 2z-l+1)(\pm 2z-l+2)\cdots(\pm 2z)}
F(l,\pm z,\pm 2z+1,1),
\end{eqnarray*}
where $F(\alpha,\beta,\gamma;t)$ denotes the Gauss hypergeometric  
function.
Using Gauss' formula for $F(\alpha,\beta,\gamma;1)$ gives
$$
\sum_{n=l}^{\infty} A_{n,l}^{(\pm)} =\pm \frac{(-1)^l}{z}.
$$
Hence we have
$$
\sum_{n=1}^{\infty} a_n=
\frac{1}{z}\sum_{l=1}^{\infty} (-1)^{l}
\left( \frac{1}{x+z-l} -\frac{1}{x-z-l} \right).
$$
Expanding the right-hand side in power series of $x$ and $z$ and taking  
the
coefficient of $x^{k-2s}z^{2s-2}$ (cf. (\ref{phi0})) gives
$$
2 {{k-1}\choose{2s-1}} \sum_{l=1}^\infty \frac{(-1)^{l-1}}{l^k},
$$
and now using the relation
$\displaystyle
\sum_{l=1}^\infty \frac{(-1)^{l-1}}{l^k}=(1-2^{1-k})\zeta(k)
$
yields equation (\ref{maintheo}).


\vspace{2mm}

\noindent{\bf\large Appendix}

\vspace{3mm}
\noindent The relation given in Theorem 1 can be interpreted as an  
equality concerning an integral
which contains the Gauss hypergeometric function.

\begin{theo} Under suitable conditions for the parameters $x$ and $z$  
that guarantee existence of
both members, the following equality holds:
\begin{eqnarray}\label{integ}
\frac{1}{1-x}\int_0^1 (1-t)^{z-x} F(1-x+z,1+z,2-x;t)\,  
dt\quad\quad\quad \quad\quad \nonumber\\
\quad\quad\quad \quad\quad =\displaystyle  
\frac{1}{z}\sum_{l=1}^{\infty} (-1)^{l}
\left( \frac{1}{x+z-l} -\frac{1}{x-z-l} \right).
\end{eqnarray}
\end{theo}

\vspace{2mm}
\noindent{\bf Proof}  \, We set
\begin{eqnarray*}
\phi_1(t)&=&t^{x+z} F(x+z,z,2z+1;t), \\
\phi_2(t)&=&t^{x-z} F(x-z,-z,-2z+1;t).
\end{eqnarray*}
Then $(\phi_1, \phi_2)$ is a system of fundamental solutions of the  
homogeneous
equation of (\ref{dfe}). The unique holomorphic solution $\Phi_0$ of  
(\ref{dfe}) is
constructed in the form
$$
\Phi_0 = u_1 \phi_1 + u_2 \phi_2,
$$
  where $u_1$ and $u_2$ are defined as follows:
\begin{eqnarray*}
  u_1(t)&=&\frac{1}{2 z}\int_0^t s^{-x-z}(1-s)^{x-1} F(x-z,-z,-2z+1;s)\,  
ds, \\
  u_2(t)&=&-\frac{1}{2z}\int_0^t s^{-x+z}(1-s)^{x-1} F(x+z,z,2z+1;s)\,  
ds.
  \end{eqnarray*}
   The values $\phi_1(1)$ and $\phi_2(1)$ are obtained by using Gauss'  
formula.
   Hence we find that $\Phi_0(1)=u_1(1)\phi_1(1)+u_2(1)\phi_2(1)$ has  
the value
\begin{eqnarray*}
   \int_0^1 dt &(1-t)^{x-1}t^{-x+z}\Big\{\displaystyle
   \frac{\Gamma(-2z)\Gamma(1-z)}{\Gamma(1-x-z)\Gamma(1-z)}
   F(x+z,z,2z+1;t)\\
   &  
\quad\quad\displaystyle+\frac{\Gamma(2z)\Gamma(1-z)}{\Gamma(1- 
x+z)\Gamma(1+z)}
   t^{z-x}F(x-z,-z,-2z+1;t)\Big\}.
  \end{eqnarray*}
Using one of the connection formulas for the Gauss hypergeometric  
functions (e.g., (43), p.~108 in \cite{erd}) yields
$$
\Phi_0(1)=\frac{1}{1-x}\int_0^1 t^{z-x}F(1-x+z,1+z,2-x;1-t)dt,
$$
which is equal to the left-hand side of (\ref{integ}). The right-hand  
side has been already
obtained in the proof of Theorem 1. This proves Theorem 2.

\vspace{2mm}
\noindent{\bf Remark} \, The right-hand side of (\ref{integ}) can be  
written in the following
form:
$$
-\frac{1}{z}\left(\psi (1-(x+z))-\psi(1-(x-z))
-\psi \left(1-\frac{x+z}{2} \right)+\psi \left(1-\frac{x-z}{2} \right)  
\right),
$$
where $\displaystyle{\psi(t)=\frac{\Gamma'(t)}{\Gamma(t)}}$ is the  
di-gamma function.

\hfill Department of Mathematics


\hfill Kinki University

\hfill Higashi-Osaka, Osaka, 577-8502 Japan

\hfill {\it aoki@math.kindai.ac.jp}

\hfill {\it ohno@math.kindai.ac.jp}

\end{document}